\documentclass[11pt]{article}

\usepackage{amsmath, amsthm, amssymb}
\usepackage{verbatim}

\title{Two remarks on even and oddtown problems}
\date{}
\author{Benny Sudakov
\thanks{Department of Mathematics, ETH, 8092 Zurich, Switzerland.
Email: benjamin.sudakov@math.ethz.ch. Research supported in part by
SNSF grant 200021-149111.} \and Pedro Vieira \thanks{Department of
Mathematics, ETH, 8092 Zurich, Switzerland.  Email:
pedro.vieira@math.ethz.ch.}}

\theoremstyle{plain}
\newtheorem{thm}{Theorem}

\newtheorem{lemma}[thm]{Lemma}

\newtheorem*{prob*}{Problem}
\theoremstyle{definition}

\newtheorem*{dfn*}{Definition}
\newtheorem{cons}{Construction}
\newtheorem*{cons*}{Construction}

\newcommand{\NN}{\mathbb{N}}

\begin{document}
\maketitle
\begin{abstract}
A family $\mathcal A$ of subsets of an $n$-element set is called an
\textit{eventown} (resp. \textit{oddtown}) if all its sets have even
(resp. odd) size and
all pairwise intersections have even size. Using tools from linear
algebra, it was shown by Berlekamp and Graver  that the maximum size
of an eventown is $2^{\left\lfloor n/2\right\rfloor}$. On the other
hand (somewhat surprisingly), it was proven by Berlekamp, that
oddtowns have size at most $n$. Over the last four decades, many
extensions of this even/oddtown problem have been studied.
In this paper we present new results on two such extensions. First,
extending a result of Vu, we show that a $k$-wise eventown (i.e.,
intersections of $k$ sets are even) has for $k \geq 3$ a unique
extremal configuration and obtain a stability result for this problem.
Next we improve some known bounds for the defect version of an
$\ell$-oddtown problem. In this problem we consider sets of size $\not\equiv 0 \pmod \ell$ where $\ell$ is a prime number $\ell$ (not necessarily $2$) and allow a few pairwise intersections to also have size $\not\equiv 0 \pmod \ell$.
\end{abstract}

\section{Introduction}
Let $\mathcal A=\{A_1,\ldots, A_m\}$ be a family of subsets of $[n]:=\{1,2,\ldots, n\}$. We say that $\mathcal A$ is an \textit{eventown} (resp. \textit{oddtown}) if all its sets have even (resp. odd) size and
\begin{align*}
|A_i\cap A_j| \;\text{ is even}&\text{ for }1\le i< j\le m
\end{align*}
Answering a question of Erd\H{o}s, Berlekamp~\cite{B69} and Graver~\cite{G75} showed independently that the maximum size of an eventown is $2^{\left\lfloor n/2\right\rfloor}$. Somewhat surprisingly, the answer changes drastically when one considers oddtowns. Indeed, Berlekamp~\cite{B69} proved that oddtowns have size at most $n$, which is easily seen to be best possible. The proofs of these two results relied on a technique known as the \textit{linear algebra bound} method, which has been widely used to tackle problems in Extremal Combinatorics ever since.

Over the last decades, many extensions of this even/oddtown problem have been studied. A natural extension is to consider the problem modulo $\ell \ge 2$. We say that $\mathcal A$ is a $\ell$-\textit{eventown} (resp. $\ell$-\textit{oddtown}) if all its sets have size $\equiv 0 \; (\bmod\; \ell)$ (resp. $\not \equiv 0 \; (\bmod\; \ell)$) and
\begin{align}
\label{eq:leventown}
|A_i\cap A_j| \equiv 0\; (\bmod\; \ell)\, \text{ for }1\le i< j\le m.
\end{align}
The problem of estimating the maximum possible size of an $\ell$-oddtown is nowadays fairly well understood. One can modify Berlekamp's proof for oddtowns slightly to show that if $\ell$ is a prime number then an $\ell$-oddtown has size at most $n$. With a bit of effort one can prove that the same still holds when $\ell$ is a prime power and that a weaker bound of $m\le c(\ell)n$ holds in general, where $c(\ell)$ is a constant depending on $\ell$. It remains an open problem whether one can take $c(\ell) = 1$ when $\ell$ is a composite number. For further details and related problems see the excellent monograph~\cite{BF92} of Babai and Frankl.

For $\ell$-eventowns a bit less is known. A natural lower bound construction for the maximum size of an $\ell$-eventown is $2^{\lfloor n/\ell\rfloor}$. This arises from considering $\lfloor n/\ell\rfloor$ disjoint subsets $B_1,\ldots, B_{\lfloor n/\ell\rfloor}$ of $[n]$ of size $\ell$ and taking $\mathcal A=\left\{\bigcup_{i\in S} B_i: \;S\subseteq \left[\lfloor n/\ell\rfloor\right]\right\}$. It turns out surprisingly that for large $\ell$ there are significantly larger $\ell$-eventowns. Indeed, Frankl and Odlyzko~\cite{FO83} found a nice construction of $\ell$-eventowns of size at least $(c\ell)^{\lfloor n/(4\ell)\rfloor}$, where $c>0$ is an absolute constant. Their construction relies on a clever use of Hadamard matrices. In addition, they showed that any $\ell$-eventown has size at most $2^{O\left(\log \ell/\ell\right)n}$ as $n\rightarrow \infty$. These two results combined certify that the maximum possible size of an $\ell$-eventown is of order $2^{\Theta(\log \ell/\ell)n}$ as $n\rightarrow \infty$.

Our results will focus on two other extensions of the even/oddtown problem that have been considered in the past. The first one extends property (\ref{eq:leventown}) to multiple intersections. The second one is a defect version of the $\ell$-oddtown problem, obtained by relaxing condition (\ref{eq:leventown}). We shall discuss these two extensions as well as our results in the next two subsections.

\subsection{Multiple intersections}
We say that $\mathcal A=\{A_1,\ldots, A_m\}$ is a \textit{$k$-wise $\ell$-eventown} if
\begin{align}
\label{eq:kwiseleventown}
\left|\bigcap_{i\in S}A_{i}\right|\equiv 0\; (\bmod\; \ell)\text{ for every non-empty }S\subseteq [m] \text{ of size } |S| = k.
\end{align}
For simplicity, we refer to a $k$-wise $2$-eventown simply as a $k$-wise eventown. We remark that a $2$-wise eventown is not the same as an eventown, since in the former we do not require that the sets themselves have even size.

The problem of maximizing the size of $k$-wise eventowns is nowadays well understood. For $k=1$, a $k$-wise eventown $\mathcal A$ is just a family of even-sized sets. Thus, $|\mathcal A| \le \sum_{i = 0}^{\lfloor n/2\rfloor}\binom{n}{2i}=2^{n-1}$, a bound which is attained by taking $\mathcal A$ to be the family of all subsets of $[n]$ of even size. The case $k=2$ was first considered in the papers of Berlekamp~\cite{B69} and Graver~\cite{G75} who showed that the maximum size of a $2$-wise eventown is $n+1$ if $n\le 5$, $2^{\lfloor n/2\rfloor}$ if $n\ge 6$ is even and $2^{\lfloor n/2\rfloor}+1$ if $n\ge 7$ is odd. Later, Vu~\cite{V97} addressed the general case:
\begin{thm}[Vu \cite{V97}]
\label{thm:vukwiseeventown}
There is a constant $c>0$ such that for any $k\ge 2$ the maximum size of a $k$-wise eventown in a universe of size $n\ge c\log_2 k$ is $2^{\lfloor n/2 \rfloor}$ if $n$ is even and $2^{\lfloor n/2\rfloor}+k-1$ if $n$ is odd.
\end{thm}
In Extremal Combinatorics, given an extremal result like Theorem~\ref{thm:vukwiseeventown}, it is common to ask what possible extremal configurations exist. In many problems, one can classify all the extremal configurations or at least describe some structural properties of these. When there is a unique extremal configuration, it is often the case that a \textit{stability} result holds. This means that one can give a precise structural description not just of the extremal configuration but also of nearly extremal configurations. 

Given Theorem~\ref{thm:vukwiseeventown}, it is therefore natural to investigate what $k$-wise eventowns of maximum possible size look like, and whether a stability version of Theorem~\ref{thm:vukwiseeventown} exists. The next construction provides $k$-wise eventowns with the sizes indicated in Theorem~\ref{thm:vukwiseeventown}, for any $k\ge 2$ and $n\ge 2\lceil\log_{2}(k-1)\rceil$.
\begin{cons}
\label{cons:maxkwiseeventown}
(i) Let $B_1,\ldots, B_{\lfloor n/2\rfloor}$ be $\lfloor n/2\rfloor$ disjoint subsets of $[n]$ of size $2$. The family $\mathcal A=\left\{\bigcup_{i\in S} B_i: \;S\subseteq \left[\lfloor n/2\rfloor\right]\right\}$ is a $k$-wise eventown of size $2^{\lfloor n/2\rfloor}$ for every $k\in \NN$.

(ii) If $n$ is odd, let $B_1,\ldots, B_{\lfloor n/2\rfloor}$ and $\mathcal A$ be as in (i). Let $i\in[n]$ be the unique element not covered by the sets $B_1,\ldots, B_{\lfloor n/2\rfloor}$ and let $C_1,\ldots, C_{k-1}$ be any $k-1$ distinct sets in $\mathcal A$ (for this we need that $n\ge 2\lceil\log_{2}(k-1)\rceil$). If we add to $\mathcal A$ the $k-1$ sets $C_1\cup\{i\},\ldots, C_{k-1}\cup\{i\}$ then the resulting family is a $k$-wise eventown of size $2^{\lfloor n/2\rfloor}+k-1$.
\end{cons}
For $k=2$, the families considered in Construction~\ref{cons:maxkwiseeventown} are by no means the only examples of $2$-wise eventowns of maximum size. For example, for $n$ even, one can show that for any $2$-wise eventown $\mathcal A$ with even-sized sets, there exists a $2$-wise eventown $\mathcal B$ containing $\mathcal A$ of size $2^{\lfloor n/2\rfloor}$ (see, e.g., Ex. 1.1.10 of Babai-Frankl~\cite{BF92}). This allows one to produce many highly non-isomorphic $2$-wise eventowns of maximum possible size, by starting with very different looking small $2$-wise eventowns $\mathcal A$ with even-sized sets and then extending them to $2$-wise eventowns of maximum possible size.  Given this phenomena, it is natural to ask what happens for $k\ge 3$. We prove that in this case the extremal construction of a $k$-wise eventown is unique. Moreover, a stability result holds.
\begin{thm}
\label{thm:kwiseeventown}
Let $\mathcal A$ be a $k$-wise eventown on $[n]$ for some $k\ge 3$. If $|\mathcal A|> \frac{3}{4}2^{\lfloor n/2\rfloor}+(k-1)n$ and $n\ge 2\lceil \log_{2}(k-1)\rceil + 4$ then $\mathcal A$ is a subfamily of a family in Construction~\ref{cons:maxkwiseeventown}.
\end{thm}
In order to establish Theorem~\ref{thm:kwiseeventown} it will be convenient for us to consider a strengthening of (\ref{eq:kwiseleventown}). We say that $\mathcal A$ is a \textit{strong $k$-wise $\ell$-eventown} if it is a $k'$-wise $\ell$-eventown for every $k'\in \{1,2,\ldots, k\}$. The problem of estimating the maximum size of a strong $k$-wise eventown is a simple one. For $k=1$, a strong $k$-wise eventown is the same as a $k$-wise eventown and so, as mentioned earlier, its maximum possible size is $2^{n-1}$. For $k\ge 2$, a strong $k$-wise eventown is also an eventown and thus has size at most $2^{\lfloor n/2\rfloor}$. Construction~\ref{cons:maxkwiseeventown} (i) certifies that strong $k$-wise eventowns of this size exist for every $k$. As was the case with $2$-wise eventowns, there are many highly non-isomorphic strong $2$-wise eventowns of size $2^{\lfloor n/2\rfloor}$. However, as our next result shows, for $k\ge 3$ the families in Construction~\ref{cons:maxkwiseeventown} (i) are the only strong $k$-wise eventowns of size $2^{\lfloor n/2\rfloor}$ and, furthermore, a stability result holds.

\begin{thm}
\label{thm:strongkwiseeventown}
If $\mathcal A$ is a $k$-wise eventown in $[n]$ for every $k\in \NN$, then there exist disjoint even-sized subsets $B_1,\ldots, B_s$ of $[n]$ such that $\mathcal A\subseteq \{\bigcup_{i\in S}B_i:S\subseteq [s]\}$. Furthermore, for $k\ge 2$, if $\mathcal A$ is a strong $k$-wise eventown in $[n]$ but not a $(k+1)$-wise eventown then $|\mathcal A|\le 2^{\left\lfloor n/2\right\rfloor-\left(2^{k}-k-2\right)}$.
\end{thm}
We remark that strong $k$-wise eventowns which are not $(k+1)$-wise eventowns only exist for $n\ge 2^{k+1}-1$. Moreover, the upper bound in Theorem~\ref{thm:strongkwiseeventown} is best possible as there exist strong $k$-wise eventowns of size $2^{\lfloor n/2\rfloor -(2^k-k-2)}$ which are not $(k+1)$-wise eventowns for any $n\ge 2^{k+1}-1$. We discuss this in Section~\ref{section:kwiseeventowns} after proving Theorem~\ref{thm:strongkwiseeventown}.

Far less is known about the maximum possible size of (strong) $k$-wise $\ell$-eventowns when $\ell > 2$. We address this problem in Section~\ref{section:concludingremarks}.

\subsection{Defect version for $\ell$-oddtowns}
We say that $\mathcal A=\{A_1,\ldots, A_m\}$ is a \textit{$d$-defect $\ell$-oddtown} if for every $i\in [m]$ we have $|A_i|\not \equiv 0\; (\bmod\; \ell)$ and there are at most $d$ indices $j\in [m]\setminus \{i\}$ such that $|A_i\cap A_j|\not \equiv 0\;(\bmod \;\ell)$. Note that a $0$-defect $\ell$-oddtown is the same as an $\ell$-oddtown. For simplicity, we refer to a $d$-defect $2$-oddtown simply as a $d$-defect oddtown. Vu~\cite{V99} considered the problem of maximizing the size of a $d$-defect oddtown, solving it almost completely. His results imply the following:
\begin{thm}[Vu~\cite{V99}]
\label{thm:vu}The maximum size of a $d$-defect oddtown in $[n]$ is $(d+1)(n-2\lceil \log_{2}(d+1)\rceil)$, 
for any $d\ge 0$ and $n\ge d/8$.
\end{thm}
For $\ell>2$, Vu observed that the maximum size of a $d$-defect $\ell$-oddtown is at most $(d+1)n$ if $\ell$ is a prime number and at least $(d+1)(n-\ell\lceil \log_{2}(d+1)\rceil)$ for every $\ell$. Our next result improves Vu's upper bound of $(d+1)n$ on the maximum size of a $d$-defect $\ell$-odtown, when $\ell >2$ is a prime number.

\begin{thm}
\label{thm:ddefectloddtown}
Let $\ell$ be a prime number and suppose $\mathcal A$ is a $d$-defect $\ell$-oddtown in the universe $[n]$. There is a constant $C>0$ such that if $n\ge Cd\log d$ then $|\mathcal A|\le (d+1)\left(n-2\left(\lceil \log_{2}(d+2)\rceil-1\right)\right)$.
\end{thm}
For $d=1$ we can show that this upper bound is essentially best possible:
\begin{thm}
\label{thm:1defectloddtown}
Let $\ell$ be a prime number. If $\mathcal A$ is a $1$-defect $\ell$-oddtown in $[n]$ then $|\mathcal A|\le \max\{n,2n-4\}$. Moreover, there exist $1$-defect $\ell$-oddtowns of size $2n-4$ for infinitely many values of $n$.
\end{thm}
It turns out that Vu's lower bound of $(d+1)(n-\ell\lceil \log_{2}(d+1)\rceil)$ can also be improved for some values of $d$ and $\ell$. We discuss this briefly in the last section of the paper.

\medskip
\noindent
\textbf{Organization of the paper:} In Section~\ref{section:auxiliaryresults} we introduce some auxiliary lemmas which we need in the proofs of our results. In Section~\ref{section:kwiseeventowns} we present the proofs of Theorems~\ref{eq:kwiseleventown} and \ref{thm:strongkwiseeventown}. In Section~\ref{section:ddefectloddtowns} we prove Theorems~\ref{thm:ddefectloddtown} and \ref{thm:1defectloddtown}. Finally, in Section~\ref{section:concludingremarks} we discuss further extensions of the problems considered as well as related open problems.

\section{Auxiliary results}
\label{section:auxiliaryresults}

The following lemma (see, e.g. Ex. 1.1.8 of \cite{BF92}) will be useful for us in the proof of Theorem~\ref{thm:kwiseeventown}.
\begin{lemma}[Skew Oddtown Theorem]
\label{lemma:skewoddtown}
Suppose $R_1,\ldots, R_m$ and $B_1,\ldots, B_m$ are subsets of $[n]$ such that the following conditions hold:
\begin{enumerate}
\item[$(a)$] $|R_i\cap B_i|\not\equiv 0 \pmod 2$ for every $i\in [m]$;
\item[$(b)$] $|R_i\cap B_j|\equiv 0 \pmod 2$ for $1\le i< j\le m$.
\end{enumerate}
Then $m\le n$.
\end{lemma}

For any graph $G$ we denote by $\chi(G)$ and $\Delta(G)$ the chromatic number and maximum degree of $G$, respectively. Recall that for any graph $G$ one has $\chi(G)\le \Delta(G)+1$ (see, e.g., \cite{D10}). In the proof of Theorem~\ref{thm:ddefectloddtown} we will be interested in the cases in which equality holds. For that matter we make use of Brooks' Theorem~\cite{B41}.
\begin{thm}[Brooks' Theorem]
\label{thm:Brooks}
For any graph $G$, we have $\chi(G)\le \Delta(G)$ unless $G$ contains a copy of $K_{\Delta(G)+1}$ or $\Delta(G)=2$ and $G$ contains a cycle of odd length.
\end{thm}

The next auxiliary lemmas use basic linear algebra. All the vector spaces considered will be over the field $\mathbb{F}_{\ell}$ where $\ell$ is a prime number and the dot product considered will always refer to the standard inner product such that $(x_1,\ldots,x_n)\cdot (y_1,\ldots, y_n) = \sum_{i=}^{n}x_iy_i$ for $(x_1,\ldots, x_n),(y_1,\ldots,y_n)\in \mathbb F_{\ell}^{n}$. We will say that a subspace $U$ of $\mathbb F_{\ell}^{n}$ is \textit{non-degenerate} if the dot product in $U$ is a non-degenerate bilinear form, meaning that for any non-zero vector $u\in U$ there exists $v\in U$ such that $u\cdot v \neq 0$. The next well-known lemma follows from Proposition 1.2 of Chapter XV of \cite{L02}.
\begin{lemma}
\label{lemma:linearalgebra}
Let $V$ be a non-degenerate subspace of $\mathbb{F}_{\ell}^{n}$ and $U$ a subspace of $V$. Denote by $U^{\perp}$ the orthogonal complement of $U$ in $V$ with respect to the dot product. Then:
\begin{enumerate}
\item[$(a)$] $\dim U+\dim U^{\perp}=\dim V$.
\item[$(b)$] If $U$ is non-degenerate then $U^{\perp}$ is also non-degenerate.
\end{enumerate}
\end{lemma}

Note that any $d$ linearly independent vectors in $\mathbb F_{\ell}^{n}$ span a subspace of size $\ell^{d}$. Therefore, given $t$ distinct vectors $v_1,\ldots, v_t$ in $\mathbb F_{\ell}^{n}$ one can always find $\lceil\log_{\ell}t\rceil$ of them which are linearly independent (e.g. take a basis of the subspace spanned by $v_1,\ldots, v_t$ consisting of vectors from this set). This is best possible in general but it can be improved under certain conditions on these vectors. A good example of this, is the following theorem of Odlyzko \cite{O81} which will be useful for us.
\begin{thm}
\label{thm:(0,1)-vectors}
Let $\ell$ be a prime number and $n$ a natural number. Given $t$ distinct $\{0,1\}$-vectors in $\mathbb F_{\ell}^{n}$ one can find at least $\lceil\log_{2}t\rceil$ of them which are linearly independent.
\end{thm}

In the proof of Theorem~\ref{thm:ddefectloddtown} we will make use of the following lemma of this type.

\begin{lemma}
\label{lemma:dimension}
Suppose $b_1,\ldots,b_{t}$ are distinct $\{0,1\}$-vectors in a non-degenerate subspace $W$ of $\mathbb{F}_{\ell}^{n}$ such that $(b_1\cdot b_1)(b_i\cdot b_j)=(b_1\cdot b_i)(b_1\cdot b_j)\neq 0$ for every $i,j\in [t]$. Then $\dim W\ge 2\lceil \log_{2}(t+1)\rceil-1$.
\end{lemma}

\begin{proof}
For each $i\in [t]$ define $c_i:=(b_1\cdot b_1)b_{i}-(b_1\cdot b_i)b_{1}$. Let $B$ and $C$ be the linear subspaces generated by $b_1,\ldots, b_{t}$ and $c_1,\ldots, c_{t}$, respectively, and let $C^{\perp}$ denote the orthogonal complement of $C$ in $W$ with respect to the dot product. Note that
\[c_{i}\cdot b_{j}=(b_1\cdot b_1)(b_i\cdot b_j)-(b_1\cdot b_i)(b_1\cdot b_j)=0\]
for every $i,j\in [t]$ and so it follows that $C\subseteq B\subseteq C^{\perp}$. Moreover, we know that $b_1\notin C$ since $b_1\cdot b_1\neq 0$ and so $\dim C\le \dim B -1$. In addition, by the definition of the vectors $c_1,\ldots, c_{t}$ it follows that $B=C+\text{span}(b_1)$ and so $\dim C\ge \dim B-1$. We conclude then that $\dim C=\dim B-1$.

By (a) of Lemma~\ref{lemma:linearalgebra} we have $\dim C+\dim C^{\perp}=\dim W$ and so we get:
\[\dim W \ge \dim B+\dim C=2\dim B-1\]
Finally, since $b_1,\ldots, b_{t}$ and the $0$-vector are $t+1$ distinct $\{0,1\}$-vectors (because $b_i\cdot b_i\neq 0$) it follows from Theorem~\ref{thm:(0,1)-vectors} that $\dim B\ge \lceil\log_{2}(t+1)\rceil$.
\end{proof}

\noindent
\textbf{Remark:} For $\ell = 2$, since all the vectors in $\mathbb{F}_{2}^{n}$ are $\{0,1\}$-vectors, one can apply Theorem~\ref{thm:(0,1)-vectors} to the vectors in $C$ to get the stronger bound $\dim W\ge 2\lceil\log_{2}t\rceil+1$. We believe that one should be able to get the same bound for any prime $\ell$.

\section{$k$-wise eventowns}
\label{section:kwiseeventowns}

In this section we present the proofs of Theorems~\ref{thm:kwiseeventown} and \ref{thm:strongkwiseeventown}. The main ingredient in the proof of Theorem~\ref{thm:kwiseeventown} is the structure of large strong $k$-wise eventowns obtained from Theorem~\ref{thm:strongkwiseeventown}. Therefore, we start with the proof of the latter and later use it to deduce the proof of the former.

\subsection{Proof of Theorem~\ref{thm:strongkwiseeventown}}
\label{subsection:strongkwiseeventown}
In the next lemma, we prove the first half of the statement in Theorem~\ref{thm:strongkwiseeventown}, characterizing the families which are $k$-wise eventowns for every $k\in \NN$.

\begin{lemma}
\label{lemma:blockfamily}
If $\mathcal A$ is a $k$-wise eventown for every $k\in \NN$, then there exist disjoint even-sized subsets $B_1,\ldots, B_s$ of $[n]$ such that $\mathcal A\subseteq \{\bigcup_{i\in S}B_i:S\subseteq [s]\}$.
\end{lemma}
\begin{proof}
Suppose $\mathcal A=\{A_1,\ldots, A_m\}$ is a $k$-wise eventown for every $k\in \NN$. Define for each $i\in [m]$ the sets $A^{0}_{i}:=A_i$ and $A^{1}_{i}:=[n]\setminus A_i$. Set $\mathcal T =\{0,1\}^{m}\setminus\{(1,1,\ldots,1)\}$ and given a tuple $t=(t_i)_{i\in [m]}\in \mathcal T$ let $B_t:=\bigcap_{i\in [m]}A^{t_i}_i$. To prove Lemma~\ref{lemma:blockfamily} it suffices to show that the sets $\{B_t:t\in \mathcal T\}$ satisfy:
\begin{enumerate}
\item[(a)] for every $i\in [m]$ there exists a set $T_i\subseteq \mathcal T$ such that $A_i=\cup_{t\in T_i}B_t$
\item[(b)] for any $t,t'\in \mathcal T$, if $t\neq t'$ then $B_t\cap B_{t'}=\emptyset$.
\item[(c)] $|B_t|$ is even for every $t\in \mathcal T$.
\end{enumerate}
We start by showing that (a) holds. Given $i\in [m]$ let $T_i=\{t\in \mathcal T: t_i=0\}$. Note that for any $t\in T_i$ we have $B_t=\bigcap_{j\in [m]}A^{t_j}_j\subseteq A_i$ since the term $A_i^{t_i}=A_i$ appears in this intersection. Thus, it follows that $\bigcup_{t\in T_i}B_t\subseteq A_i$. Now, note that for each $a\in A_i$ there exists $t\in T_i$ such that $a\in B_t$. Indeed, just consider $t_j=0$ if $a\in A_j$ and $t_j=1$ otherwise. Thus, it follows also that $A_i\subseteq \bigcup_{t\in T_i}B_t$.

Next, we show that (b) holds. Suppose $t\neq t'$ and let $i\in [m]$ be such that $t_i\neq t'_i$. Then $B_t\subseteq A_i^{t_i}$ and $B_{t'}\subseteq A_i^{t'_i}$. Since $t_i\neq t'_i$ it follows that $A_i^{t_i}\cap A_i^{t'_i}=\emptyset$ and so $B_t\cap B_{t'}=\emptyset$.

Finally, we show that $(c)$ holds. Given $t\in \mathcal T$ we have:
\begin{align*}|B_t|=&\left|\left(\bigcap_{i\in [m], t_i=0}A_i\right)\cap\left(\bigcap_{i\in [m],t_i=1}[n]\setminus A_i\right)\right|\\
=&\left|\left(\bigcap_{i\in [m], t_i=0}A_i\right)\setminus\left(\bigcup_{i\in [m],t_i=1}A_i\right)\right|\\
=& \left|\left(\bigcap_{i\in [m], t_i=0}A_i\right)\right|-\left|\left(\bigcap_{i\in [m], t_i=0}A_i\right)\cap \left(\bigcup_{i\in [m],t_i=1}A_i\right)\right|.
\end{align*}
The first term is the intersection of a positive number of sets in $\mathcal A$ (since $t\neq (1,1,\ldots,1)$) and thus has even size since $\mathcal A$ is a $k$-wise eventown for every $k\in \NN$. Moreover, the second term can be written, by the inclusion-exclusion principle, as a sum of signed intersection sizes of sets in $\mathcal A$. Thus, the second term is also even, implying that $|B_t|$ is even.
\end{proof}

\medskip

For the second half of the statement of Theorem~\ref{thm:strongkwiseeventown} we will use basic linear algebra techniques. Given a set $A\subseteq [n]$ let $v_A\in\mathbb{F}_{2}^{n}$ denote its $\{0,1\}$-characteristic vector. We consider the following two correspondences between families $\mathcal A\subseteq 2^{[n]}$ and linear subspaces $V\subseteq\mathbb{F}_{2}^{n}$:
\[\mathcal A\mapsto V_{\mathcal A}:=\text{span}\{v_{A}:A\in\mathcal A\}\;\;\text{and}\;\;V\mapsto \mathcal A_V:=\{A\subseteq [n]: v_{A}\in V\}\]
Given $\mathcal A\subseteq 2^{[n]}$, we define $\overline{\mathcal A}:=\mathcal A_{V_{\mathcal A}}$ which we call the \textit{linear closure of $\mathcal A$}. Note that $\mathcal A\subseteq \overline{\mathcal A}$, but equality does not necessarily hold. As the next lemma shows, an important property of linear closure is that it preserves the property of being a strong $k$-wise eventown.
\begin{lemma}
\label{lemma:linearclosure}
If $\mathcal A$ is a strong $k$-wise eventown then $\overline{\mathcal A}$ is also a strong $k$-wise eventown.
\end{lemma}

\begin{proof}
Given a set $B\subseteq [n]$ define the function $f_{B}:[n]\rightarrow \mathbb F_{2}$ such that 
\[f_{B}(i)=\left\{\begin{matrix}
1 & \text{ if } i\in B\\ 
0 & \text{ if } i\notin B
\end{matrix}\right.\]
and note that:
\begin{enumerate}
\item[(i)] for any $B\subseteq [n]$ we have $|B|= \sum_{i\in [n]}f_B(i)\pmod 2$;
\item[(ii)] for any $t$ sets $B_1,\ldots,B_t\subseteq [n]$ we have $f_{\cap_{i\in [t]}B_i}=\prod_{i\in [t]}f_{B_i}$;
\item[(iii)] if $A_1,\ldots, A_{t}, B\subseteq [n]$ are such that $v_{B}=\sum_{i\in [t]}v_{A_i}$ then $f_{B}=\sum_{i\in [t]}f_{A_{i}}$.
\end{enumerate}
Now, let $B_1,\cdots,B_k$ be any $k$ not necessarily distinct sets in $\overline{\mathcal A}$. We want to show that $\bigcap_{j\in [k]}B_j$ has even size. Since $\overline{\mathcal A}$ is the span of the vectors $\{v_A\}_{A\in \mathcal A}$, we know that for each $j\in [k]$ there are sets $A_1^{j},\ldots, A_{t_{j}}^{j}\in \mathcal A$ such that $v_{B_j}=\sum_{i\in [t_j]}v_{A_{i}^{j}}$. Thus, by properties (i), (ii) and (iii) it follows that
\begin{align*}
\left|\bigcap_{j\in[k]} B_j\right|&= \sum_{i\in [n]}f_{\cap_{j\in[k]} B_j}(i)\\
&= \sum_{i\in [n]}\prod_{j\in [k]}f_{B_j}(i)\\
&= \sum_{i\in [n]}\prod_{j\in [k]}\sum_{h\in [t_j]}f_{A_h^{j}}(i)\\
&=\sum_{i\in [n]}\sum_{(h_1,\ldots, h_k)}f_{\cap_{j\in [k]}A_{h_j}^{j}}(i)\\
&= \sum_{(h_1,\ldots,h_k)}\sum_{i\in [n]}f_{\cap_{j\in [k]}A_{h_j}^{j}}(i)\\
&= \sum_{(h_1,\ldots,h_k)}\left|\cap_{j\in [k]}A_{h_j}^{j}\right|\pmod 2
\end{align*}
where the sums indexed with $(h_1,\ldots,h_k)$ run over all tuples in $[t_1]\times\ldots\times[t_k]$.
Since $\mathcal A$ is a strong $k$-wise eventown we conclude that all the terms in the last sum are even. Thus, for any $k$ not necessarily distinct sets $B_1,\ldots, B_k\in\overline{\mathcal A}$ the set $\bigcap_{j\in[k]} B_j$ has even size, i.e., $\overline{\mathcal A}$ is a strong $k$-wise eventown.
\end{proof}

With Lemma~\ref{lemma:linearclosure} we are ready to present the proof of the second half of the statement of Theorem~\ref{thm:strongkwiseeventown}.

\begin{lemma}
\label{lemma:upperbound}
For $k\ge 2$, if $\mathcal A\subseteq 2^{[n]}$ is a strong $k$-wise eventown but not a $(k+1)$-wise eventown then:
\[|\mathcal A|\le 2^{\left\lfloor n/2\right\rfloor-\left(2^{k}-k-2\right)}\]
\end{lemma}

\begin{proof}[Proof of Lemma~\ref{lemma:upperbound}]
Suppose $\mathcal A\subseteq 2^{[n]}$ is a strong $k$-wise eventown which is not a $(k+1)$-wise eventown and let $A_1,\ldots, A_{k+1}\in \mathcal A$ be such that $|A_1\cap\ldots\cap A_{k+1}|$ is odd. For each $S\subseteq [k+1]$ define the set $A_{S}:=\bigcap_{i\in S}A_i$, let $\mathcal S=\{S\subseteq [k]: 2\le |S|\le k-1\}$ and define $\mathcal B:=\{A_{S}\}_{S\in\mathcal S}$. We claim that the family $\mathcal C=\mathcal A\cup\mathcal B$ is an eventown. Indeed, this holds since:
\begin{enumerate}
\item[1)] all sets in $\mathcal A$ and pairwise intersections between sets in $\mathcal A$ have even size since $\mathcal A$  is a strong $k$-wise eventown and $k\ge 2$;
\item[2)] all sets in $\mathcal B$ have even size since they are the intersection of at most $k-1$ sets in $\mathcal A$;
\item[3)] for any $A\in \mathcal A$ and $S\in\mathcal S$ the set $A\cap A_{S}=A\cap \left(\bigcap_{i\in S}A_i\right)$ is the intersection of at most $k$ sets in $\mathcal A$, and thus has even size;
\item[4)] for any $S_1,S_2\in\mathcal S$ the set $A_{S_1}\cap A_{S_2}=\bigcap_{i\in S_1\cup S_2}A_i$ is the intersection of at most $k$ sets in $\mathcal A$, and thus has even size.
\end{enumerate}
We claim now that $\dim V_{\mathcal C}=\dim V_{\mathcal A} + \dim V_{\mathcal B}$ and that $\dim V_{\mathcal B} = |\mathcal S|=2^{k}-k-2$. If this is the case then:
\[\left|\overline{\mathcal C}\right|=2^{\dim V_{\mathcal C}}=2^{\dim V_{\mathcal A}}\cdot 2^{\dim V_{\mathcal B}}\ge |\mathcal A|\cdot 2^{2^{k}-k-2}\]
and since $\overline{\mathcal C}$ is an eventown by Lemma~\ref{lemma:linearclosure}, we conclude that
\[|\mathcal A|\le \left|\overline{\mathcal C}\right|\cdot 2^{-\left(2^{k}-k-2\right)}\le 2^{\left\lfloor n/2\right\rfloor-\left(2^{k}-k-2\right)}\]
as desired. Thus, it remains to prove the claim. For that, it suffices to prove that if there is a linear relation
\begin{align}
\label{eq:lindep2}
\sum_{A\in \mathcal A}\alpha_{A}v_{A}+\sum_{S\in \mathcal S}\beta_{S}v_{A_S}=0
\end{align}
then $\beta_{S}=0$ for any $S\in \mathcal S$. Define for each $S\in \mathcal S$ the set $S^{c}:=[k+1]\setminus S$ and note that for any $A\in \mathcal A$ and $S,T\in \mathcal S$ we have:
\begin{enumerate}
\item[(i)] $v_{A}\cdot v_{A_{T^{c}}}=|A\cap\left(\bigcap_{i\in T^{c}}A_i\right)|= 0\pmod 2$ because the latter is the intersection of at most $k$ sets in $\mathcal A$, since $|T^{c}|=k+1-|T|\le k-1$;
\item[(ii)] if $S\cup T^{c}\neq [k+1]$ then $v_{A_S}\cdot v_{A_{T^{c}}}=|\bigcap_{i\in S\cup T^{c}}A_i|= 0\pmod 2$ because the latter is the intersection of at most $k$ sets in $\mathcal A$;
\item[(iii)] $v_{A_T}\cdot v_{A_{T^{c}}}=|\bigcap_{i\in[k+1]}A_i|= 1\pmod 2$.
\end{enumerate}
Consider now a linear relation as in equation~(\ref{eq:lindep2}) and suppose that there is some set $S\in \mathcal S$ such that $\beta_{S}\neq 0$. Let $T\in \mathcal S$ be such a set of maximum possible size and note that for any $S\in \mathcal S\setminus \{T\}$ with $\beta_{S} \neq 0$ we have $T\not\subseteq S$, or equivalently $
S\cup T^{c}\neq [k+1]$. Therefore, it follows from (i), (ii) and (iii) that
\[0 = \left(\sum_{A\in \mathcal A}\alpha_{A}v_{A}+\sum_{S\in \mathcal S}\beta_{S}v_{A_S}\right) \cdot v_{A_{T^{c}}}=\sum_{A\in \mathcal A}\alpha_{A}\left(v_{A}\cdot v_{A_{T^{c}}}\right)+\sum_{S\in \mathcal S}\beta_{S}\left(v_{A_S}\cdot v_{A_{T^{c}}}\right) =\beta_{T}\]
contradicting the choice of $T$. This proves the claim.
\end{proof}

Note that Lemma~\ref{lemma:upperbound} implies that there is no strong $k$-wise eventown in $[n]$ that is not a $(k+1)$-wise eventown if $\lfloor n/2\rfloor < 2^{k}-k-2$. In fact, one actually needs that $n\ge 2^{k+1}-1$ for such families to exist. The reason for this is quite simple. If $\mathcal A$ is not a $(k+1)$-wise eventown then there exist sets $A_1,\ldots, A_{k+1}\in \mathcal A$ for which $|A_1\cap \ldots, A_{k+1}|$ is odd. Since the intersection of the sets in any proper non-empty subfamily of $\{A_1,\ldots, A_{k+1}\}$ has even size then one can use the principle of inclusion-exclusion to show that in fact $|A'_1\cap \ldots \cap A'_{k+1}|$ is odd for any choice of $A'_i \in \{A_i, [n]\setminus A_i\}$ for $i\in [k+1]$, with the exception of the choice $A'_i = [n]\setminus A_i$ for every $i\in [k+1]$ (when $n$ is odd). This implies that there are at least $2^{k+1}-1$ disjoint non-empty sets in $[n]$, implying that $n\ge 2^{k+1}-1$.

We show next that for any $n\ge 2^{k+1}-1$ there are strong $k$-wise eventowns $\mathcal A$ in $[n]$ of size $|\mathcal A|=2^{\lfloor n/2\rfloor-(2^{k}-k-2)}$ which are not $(k+1)$-wise eventowns. We start by constructing a strong $k$-wise eventown consisting of $2^{k+2}$ subsets of $[2^{k+1}]$ which is not a $(k+1)$-wise eventown.

For convenience, let us denote by $2^{[k+1]}$ the family of all subsets of the set $[k+1] = \{1,\ldots, k+1\}$ and let $f:2^{[k+1]}\rightarrow [2^{k+1}]$ be any bijection. Let $B_0=[2^{k+1}]$ and for each $i\in [k+1]$ define $B_i=\{f(S):i\in S\subseteq [k+1]\}$. Note that for any set $I\subseteq \{0,1,\ldots,k+1\}$ we have:
\[\left|\bigcap_{i\in I}B_i\right|=\left|\bigcap_{i\in I\setminus \{0\}}B_i\right|=\left|\left\{f(S):\left(I\setminus \{0\}\right)\subseteq S\subseteq [k+1]\right\}\right|=2^{k+1-\left|I\setminus \{0\}\right|}\]
and so the family $\mathcal B=\{B_0,B_1,\ldots, B_{k+1}\}$ is a strong $k$-wise eventown but not a $(k+1)$-wise eventown. Hence, by Lemma~\ref{lemma:linearclosure} it follows that $\overline{\mathcal B}$, the linear closure of $\mathcal B$, is also a strong $k$-wise eventown but not a $(k+1)$-wise eventown.

We claim now that the vectors $v_{B_0},\ldots, v_{B_{k+1}}$ are linearly independent. Indeed, this follows from the next observations:
\begin{itemize}
\item $v_{\{f(\emptyset)\}}\cdot v_{B_0} = 1$ and $v_{\{f(\emptyset)\}}\cdot v_{B_i} = 0$ for $i\in [k+1]$ since $f(\emptyset)\not\in B_i$.
\item for $i,j\in[k+1]$ we have $v_{\{f(\{i\})\}}\cdot v_{B_j} = \left\{\begin{matrix}
1 & \text{if }i=j\\ 
0 & \text{if }i\neq j
\end{matrix}\right.$.
\end{itemize}
Therefore, $\overline{\mathcal B}$ is a strong $k$-wise eventown in $[2^{k+1}]$ of size $|\overline{\mathcal B}| = 2^{\dim V_{\mathcal B}} = 2^{k+2}$ which is not a $(k+1)$-wise eventown.

Now, if $n\ge 2^{k+1}$ let $\mathcal C$ be a strong $k$-wise eventown in $[n]\setminus [2^{k+1}]$ of size $2^{\lfloor (n-2^{k+1})/2\rfloor}$ as in Construction~\ref{cons:maxkwiseeventown} (i). Since $\overline{\mathcal B}$ and $\mathcal C$ are both strong $k$-wise eventowns and $\overline{\mathcal B}$ is not a $(k+1)$-wise eventown, a moment's thought reveals that the family
\[\mathcal A = \{B\cup C: B\in \overline{\mathcal B}, C\in \mathcal C\}\]
is a strong $k$-wise eventown in $[n]$ of size
\[|\mathcal A|= 2^{k+2}\cdot 2^{\lfloor (n-2^{k+1})/2\rfloor} = 2^{\lfloor n/2\rfloor-(2^{k}-k-2)}\]
which is not a $(k+1)$-wise eventown.

When $n=2^{k+1}-1$, note that if we choose the bijection $f$ above such that $f(\emptyset) = 2^{k+1}$ then the sets $B_1,\ldots, B_{k+1}$ are subsets of $[2^{k+1}-1] = [n]$. Therefore, in a similar way as above, we can conclude that the linear closure of $\{B_1,\ldots, B_{k+1}\}$ will be a strong $k$-wise eventown in $[n]$ of size $2^{k+1} = 2^{\lfloor n/2\rfloor - (2^{k}-k-2)}$ which is not a $(k+1)$-wise eventown.

\subsection{Proof of Theorem~\ref{thm:kwiseeventown}}
We will use Theorem~\ref{thm:strongkwiseeventown} in order to prove Theorem~\ref{thm:kwiseeventown}. The reason why we can do this is because, as the next lemma shows, any $k$-wise eventown contains a large strong $k$-wise eventown.
\begin{lemma}
\label{lemma:strongsubfamilies}
If $\mathcal A$ is a $k$-wise eventown on $[n]$ then it contains a subfamily $\mathcal A'$ of size $|\mathcal A'|\ge |\mathcal A|-(k-1)n$ which is a strong $k$-wise eventown.
\end{lemma}
\begin{proof}
Set $\mathcal A_0 := \mathcal A$ and, for $i\ge 0$, as long as $\mathcal A_{i}$ is not a strong $k$-wise eventown let $A^{i}_1,\ldots,A^{i}_{k_{i}}$ be a maximal collection of less than $k$ distinct sets in $\mathcal A_{i}$ such that $|A^{i}_1\cap\ldots\cap A^{i}_{k_{i}}|\not\equiv 0\pmod 2$ and set $\mathcal A_{i+1}:=\mathcal A_{i}\setminus\{A^{i}_1,\ldots, A^{i}_{k_i}\}$. After a finite number of iterations of this procedure, say $s$ iterations, we obtain a (possibly empty) subfamily $\mathcal A'$ of $\mathcal A$ which is a strong $k$-wise eventown. Since at each step $i<s$ the family $\mathcal A_{i+1}$ is obtained from $\mathcal A_i$ by removing $k_i \le k-1$ sets, we have:
\[|\mathcal A'|\ge |\mathcal A|-(k-1)s.\]
Thus, it suffices then to show that $s\le n$. 

For each $i\in [s]$ define the sets $R_i:=A_1^{i}\cap\ldots A_{k_i}^{i}$ and $B_i:=A_1^{i}$ and note that:
\begin{enumerate}
\item[(a)] $|R_i\cap B_i|=|A_1^{i}\cap\ldots A_{k_i}^{i}|\not\equiv 0 \pmod 2$ for every $i\in [s]$;
\item[(b)] $|R_i\cap B_j|=|A_1^{i}\cap\ldots A_{k_i}^{i}\cap A_1^{j}|\equiv 0 \pmod 2$ for $i< j$ since otherwise $A_1^{i},\ldots, A_{k_i}^{i}, A_1^{j}$ would be a collection of $k_i+1$ distinct sets in $\mathcal A_i$ whose intersection has odd size, contradicting the maximality in the choice of the sets $A_1^{i},\ldots, A_{k_i}^{i}$ (note that $k_i + 1 < k$ since $\mathcal A$ is a $k$-wise eventown).
\end{enumerate}
Thus, by Lemma~\ref{lemma:skewoddtown} it follows that $s\le n$, as desired.
\end{proof}

\begin{proof}[Proof of Theorem~\ref{thm:kwiseeventown}]
By Lemma~\ref{lemma:strongsubfamilies}, the family $\mathcal A$ contains a subfamily $\mathcal A'$ of size
\begin{align}
\label{eq:sizeofsubfamily}
|\mathcal A'|\ge |\mathcal A|-(k-1)n> \frac{3}{4}2^{\lfloor n/2\rfloor}
\end{align}
which is a strong $k$-wise eventown. Thus, since $2^{\lfloor n/2\rfloor -(2^{k}-k-2)}\le \frac{1}{8}2^{\lfloor n/2\rfloor}\le\frac{3}{4}2^{\lfloor n/2\rfloor}$, it follows from Theorem~\ref{thm:strongkwiseeventown} that there are non-empty disjoint subsets $B_1,\ldots, B_s$ of $[n]$ of even size such that $\mathcal A'\subseteq \mathcal B:=\{\bigcup_{i\in S}B_i : S\subseteq [s]\}$. Note that $|\mathcal A'|\le |\mathcal B|\le 2^{s}$ and so it follows from (\ref{eq:sizeofsubfamily}) that $s\ge\lfloor n/2\rfloor$. Furthermore, since the sets $B_1,\ldots, B_s$ are non-empty disjoint subsets of $[n]$ of even size, we must have $s=\lfloor n/2\rfloor$ and $|B_i|=2$ for every $i\in [s]$.

We claim now that for any $A\in \mathcal A$ and $i\in [s]$, if $A\cap B_i\neq \emptyset$ then $B_i\subseteq A$. Suppose that this is not the case and let $A^{*}\in \mathcal A$ and $i\in [s]$ be such that $|A^{*}\cap B_i|=1$. Let $\mathcal A''=\{A\in \mathcal A': B_i\subseteq A\}$ and note that $\mathcal A'\setminus \mathcal A'' \subseteq \{\bigcup_{j\in S}B_j: S\subseteq [s]\setminus \{i\}\}$. Therefore $|\mathcal A'\setminus \mathcal A''|\le 2^{s-1}$ and so by (\ref{eq:sizeofsubfamily}):
\[|\mathcal A''|\ge |\mathcal A'|-2^{s-1}=|\mathcal A'|-\frac{1}{2}2^{s}> \frac{1}{4}2^{s} = \frac{1}{2}2^{s-1}.\] 
Thus, if we define $\mathcal S=\left\{S\subseteq [s]\setminus \{i\}:B_i\cup\left(\bigcup_{j\in S}B_j\right)\in \mathcal A''\right\}$, we see that $\mathcal S\subseteq 2^{[s]\setminus\{i\}}$ and that $|\mathcal S|=|\mathcal A''|>\frac{1}{2}\left|2^{[s]\setminus\{i\}}\right|$. Hence, there must exist two distinct disjoint sets $S_{1},S_{2}\subseteq [s]\setminus \{i\}$ such that $A_1:=B_i\cup\left(\bigcup_{j\in S_{1}}B_j\right)\in \mathcal A''$ and $A_2:=B_i\cup\left(\bigcup_{j\in S_{2}}B_j\right)\in \mathcal A''$. Since $S_1$ and $S_2$ are disjoint, this implies that $A_1\cap A_2=B_i$. Finally, let $A_3,\ldots, A_{k-1}$ be $k-3$ distinct sets in $\mathcal A''\setminus \{A_1,A_2\}$ and note that
\[|A^{*}\cap A_1\cap A_2\cap\ldots\cap A_{k-1}|=|A^{*}\cap B_i|=1\]
contradicting the fact that $\mathcal A$ is a $k$-wise eventown. Thus, we conclude that for any $A\in \mathcal A$ and $i\in [s]$ if $A\cap B_i\neq \emptyset$ then $B_i\subseteq A$.

If $n$ is even then $\bigcup_{j\in [s]}B_j=[n]$ and so it follows that $\mathcal A\subseteq \mathcal B$, the latter being a family in Construction~\ref{cons:maxkwiseeventown}. If $n$ is odd, then $\bigcup_{j\in [s]}B_j=[n]\setminus \{i\}$ for some $i\in [n]$, and thus $\mathcal A\subseteq \mathcal B\cup\{C\cup \{i\}:C\in \mathcal B\}$. Since the intersection of any number of sets of the form $\{C\cup\{i\}\}_{C\in \mathcal B}$ has odd size, and since $\mathcal A$ is a $k$-wise eventown, we conclude that there are at most $k-1$ sets $C\in \mathcal B$ such that $C\cup\{i\}\in\mathcal A$. Thus, we conclude that $\mathcal A$ is a subfamily of a family in Construction~\ref{cons:maxkwiseeventown}.
\end{proof}

\noindent
\textbf{Remark:} In the proof of Theorem~\ref{thm:kwiseeventown} we implicitly use the fact that $|\mathcal A''|\ge k-1$ when we consider $k-3$ distinct sets $A_3,\ldots, A_{k-1}$ from $\mathcal A''\setminus \{A_1,A_2\}$. This follows from the fact that $|\mathcal A''|\ge \frac{1}{4}2^{\lfloor n/2\rfloor}$ and the condition $n\ge 2\lceil \log_{2}(k-1)\rceil+4$ in the theorem statement.

\section{$d$-defect $\ell$-oddtowns}
\label{section:ddefectloddtowns}

\subsection{Proof of Theorem~\ref{thm:ddefectloddtown}}

Given a family of sets $\mathcal A=\{A_1,\ldots, A_{m}\}$ we define its $\ell$-auxiliary graph $G_{\ell}(\mathcal A)$ to be the simple graph with vertex set $\mathcal A$ where $A_iA_j$ is an edge if and only if $|A_i\cap A_j|\not\equiv 0\pmod \ell$. We will often abuse notation slightly and refer to the properties of $G_{\ell}(\mathcal A)$ as being properties of $\mathcal A$. In particular, we use $\Delta(\mathcal A)$, $\chi(\mathcal A)$ and $\alpha(\mathcal A)$ to denote the maximum degree, chromatic number and independence number of $G_{\ell}(\mathcal A)$, respectively.

Let $\mathcal A$ be a $d$-defect $\ell$-oddtown in $[n]$, where $\ell$ is a prime number. Note that $\Delta(\mathcal A)\le d$ and so, in particular, $\alpha(\mathcal A)\ge |\mathcal A|/(d+1)$. 
Moreover, observe crucially that an independent set in $G_{\ell}(\mathcal A)$ corresponds to an $\ell$-oddtown inside $\mathcal A$, which as we discussed in the introduction has size at most $n$. Hence, we conclude that $|\mathcal A|\le (d+1)n$. We now wish to improve this simple upper bound to $(d+1)(n-t)$ where $t=2\left(\lceil\log_{2}(d+2)\rceil-1\right)$. We consider the following two cases:
\begin{enumerate}
\item[(a)] $G_{\ell}(\mathcal A)$ contains at most $t$ copies of $K_{d+1}$

\item[(b)] $G_{\ell}(\mathcal A)$ contains more than $t$ copies of $K_{d+1}$
\end{enumerate}
and show that in any case we have $|\mathcal A|\le (d+1)(n-t)$.

We consider case (a) first. Let $\mathcal A'$ be a family obtained from $\mathcal A$ by removing one set from each copy of $K_{d+1}$ in $G_{\ell}(\mathcal A)$. We claim that $\alpha(\mathcal A')\ge |\mathcal A'|/\left(d+\frac{1}{2}\right)$. Indeed, note that the graph $G_{\ell}(\mathcal A')$ does not contain a copy of $K_{d+1}$. Therefore, if $d\neq 2$, it follows from Brooks' Theorem (Theorem~\ref{thm:Brooks}) that $\chi(\mathcal A')\le d$, which implies that $\alpha(\mathcal A')\ge |\mathcal A'|/d\ge |\mathcal A'|/\left(d+\frac{1}{2}\right)$. If $d=2$ then, since $\Delta(\mathcal A')\le 2$, the graph $G_{\ell}(\mathcal A')$ is a disjoint union of cycles of length at least $4$ (recall that $G_{\ell}(\mathcal A')$ is $K_3$-free) and paths. A path of length $\ell$ has an independent set of size at least $\ell/2$ and a cycle of length $\ell \ge 4$ has an independent set of size at least $2\ell/5$. Thus, for $d=2$, it follows that $\alpha(\mathcal A')\ge 2|\mathcal A'|/5=|\mathcal A'|/\left(d+\frac{1}{2}\right)$. Since an independent set in $G_{\ell}(\mathcal A')$ corresponds to an $\ell$-oddtown inside $\mathcal A'$ and since an $\ell$-oddtown in $[n]$ has at most $n$ sets, we conclude that $|\mathcal{A'}|/\left(d+\frac{1}{2}\right)\le n$ and hence:
\[|\mathcal{A}|\le t+|\mathcal A'|\le t+ \left(d+\frac{1}{2}\right)n\le (d+1)(n-t)\]
provided $n\ge Cd\log d$ for some constant $C>0$.

\medskip

We consider now case (b). Let $C_1,\ldots, C_{r}$ denote the connected components of the graph $G_{\ell}(\mathcal A)$. For each $A\in \mathcal A$, let $v_{A}$ denote its characteristic vector in $\mathbb{F}_{\ell}^{n}$ and consider the $n\times |\mathcal A|$ matrix $M$ whose column vectors are the vectors $\{v_{A}\}_{A\in \mathcal A}$, ordered according to the connected components $C_1,\ldots, C_r$. Note that the matrix $\mathcal M=M^{T}M$ is a square matrix of dimension $|\mathcal A|$ and that the entry corresponding to two sets $A,B\in \mathcal A$ in $\mathcal M$ is precisely $v_{A}\cdot v_{B}=|A\cap B|\pmod \ell$. Moreover, since  the rows and columns of $\mathcal M$ are ordered according to the connected components of $G_{\ell}(\mathcal A)$ and since $|A\cap B|=0\pmod \ell$ for $A,B\in \mathcal A$ in different connected components, it follows that $\mathcal M$ is a block diagonal matrix, with each block $\mathcal M_i$ corresponding to a connected component $C_i$. Thus, we have:
\begin{align}
\label{eq: sum of ranks}
\sum_{i=1}^{r} \text{rank}(\mathcal M_i)=\text{rank}(\mathcal M)\le \text{rank}(M)\le n.
\end{align}
Note that if $\mathcal I=\{A_1,\ldots, A_{|\mathcal I|}\}$ is an independent set in $C_i$ then $v_{A_{j}}\cdot v_{A_{j'}}=|A_{j}\cap A_{j'}|\neq 0\pmod \ell$ if and only if $j=j'$, implying that the submatrix of $\mathcal M_i$ whose rows and columns correspond to the sets in $\mathcal I$ has full rank $|\mathcal I|$. Thus, since $\Delta(\mathcal A)\le d$ it follows that for each $i\in [r]$:
\begin{align}
\label{eq: bound on rank}
\text{rank}(\mathcal M_i)\ge \alpha(C_i)\ge |C_i|/(d+1).
\end{align}

We claim now that there is at least one component $C_i$ which is a copy of $K_{d+1}$ such that $\text{rank}(\mathcal M_i)=1$, or else $|\mathcal A|<(d+1)(n-t)$. Indeed, since we are looking at case (b), we know that more than $t$ components of $G_{\ell}(\mathcal A)$ are copies of $K_{d+1}$. Moreover, if all the corresponding blocks have rank at least $2$ then there are more than $t$ values of $i\in [r]$ for which inequality (\ref{eq: bound on rank}) can be improved to $\text{rank}(\mathcal M_i)\ge 1+|C_i|/(d+1)$. Thus, in that case it follows from (\ref{eq: sum of ranks}) that:
\[n\ge \sum_{i=1}^{r}\text{rank}(\mathcal M_i)> t+\sum_{i=1}^{r}|C_i|/(d+1)=t+|\mathcal A|/(d+1) \Rightarrow |\mathcal A|<(d+1)(n-t)\]

Thus, we may assume that there is one connected component $C_{i^{*}}$ of $G_{\ell}(\mathcal A)$ which is a copy of $K_{d+1}$ and whose corresponding block matrix $\mathcal M_{i^{*}}$ in $\mathcal M$ has rank $1$. Note that this implies that any two rows/columns in $\mathcal M_{i^{*}}$ are multiples of one another. Let $B_1,\ldots, B_{d+1}$ be the sets in $\mathcal A$ corresponding to such a connected component. Note that since $b_i\cdot b_i = |B_i| \neq 0 \pmod \ell$  for any $i\in [d+1]$ and since the rows of $\mathcal M_{i^{*}}$ are multiples of one another, it follows that $b_i\cdot b_j\neq 0 \pmod \ell$ for any $i,j\in [d+1]$ and that $(b_1\cdot b_1)(b_i\cdot b_j) = (b_1\cdot b_i)(b_1\cdot b_j)$.

Now, let $\mathcal A'$ denote the family $\mathcal A\setminus \{B_1,\ldots, B_{d+1}\}$ and let $A_1,\ldots, A_{s}$ be sets corresponding to an independent set of maximum size in $G_{\ell}(\mathcal A')$. Since $\Delta(\mathcal A')\le d$ it follows that
\begin{align}
\label{eq:bound}
s=\alpha(\mathcal A')\ge \frac{|\mathcal A'|}{d+1} =\frac{|\mathcal A|}{d+1}-1.
\end{align}
Let $a_1,\ldots, a_{s}$ and $b_1,\ldots,b_{d+1}$ be the characteristic vectors in $\mathbb{F}_{\ell}^{n}$ of $A_1,\ldots, A_{s}$ and $B_1,\ldots, B_{d+1}$, respectively. Because of the choices of these sets, it follows that:
\begin{enumerate}
\item[(i)] For every $i,j\in [s]$: $a_i\cdot a_j\neq 0$ if and only if $i=j$. 
\item[(ii)] For every $i,j\in [d+1]$: $(b_1\cdot b_1)(b_{i}\cdot b_j)=(b_1\cdot b_i)(b_1\cdot b_j)\neq 0$.
\item[(iii)] For every $i\in [s]$ and $j\in [d+1]$: $a_i\cdot b_j=0$.
\end{enumerate}
Denoting by $U$ the space generated by $a_1,\ldots, a_{s}$, it follows from (i) that $U$ is a non-degenerate subspace of $\mathbb{F}_{\ell}^{n}$ (see Section~\ref{section:auxiliaryresults} for the definition) and that $\dim U=s$. Furthermore, by Lemma~\ref{lemma:linearalgebra} we know that $U^{\perp}$ is non-degenerate and that
\begin{align}
\label{eq:dimension}
s+\dim U^{\perp} = \dim U + \dim U^{\perp}=n
\end{align}
Since the vectors $b_1,\ldots, b_{d+1}$ are distinct $\{0,1\}$-vectors satisfying (ii) and are in $U^{\perp}$ by (iii), we obtain by Lemma~\ref{lemma:dimension} that
\begin{align}
\label{eq:lowerbounddimension}
\dim U^{\perp}\ge 2\lceil\log_{2}(d+2)\rceil-1 = t+1.
\end{align}
Finally, putting (\ref{eq:bound}), (\ref{eq:dimension}) and (\ref{eq:lowerbounddimension}) together, we conclude that
\[\left(\frac{|\mathcal A|}{d+1}-1\right)+(t+1)\le n\; \Leftrightarrow\; |\mathcal A|\le (d+1)(n-t)\]
as claimed. This finishes the proof of Theorem~\ref{thm:ddefectloddtown}.

\subsection{Proof of Theorem~\ref{thm:1defectloddtown}}
We start by giving constructions of $1$-defect $\ell$-oddtowns of size $2n-4$ for infinitely many values of $n$, when $\ell$ is a prime number. Our constructions rely on the use of Hadamard matrices. A Hadamard matrix of order $n$ is an $n \times n$ matrix whose entries are either $+1$ or $-1$ and whose rows are mutually orthogonal. A necessary condition for a Hadamard matrix of order $n > 2$ to exist is that $n$ is divisible by $4$. The most important open question in the theory of Hadamard matrices, known as the Hadamard conjecture, is whether this condition is also sufficient. For more on Hadamard matrices see e.g. \cite{H07}.

Suppose a Hadamard matrix $H$ of order $n-1$ exists. We may assume the last column has every entry equal to $1$, by multiplying some rows by $-1$ if necessary. For $j\in [n-2]$ define sets $A_j, B_j\subseteq [n-1]$ by taking $i \in A_j$ if and only if $H_{i,j} = 1$ and setting $B_j=[n-1]\setminus A_j$. The fact that $H$ is a Hadamard matrix of order $n-1$ with the last column being the all-$1$ vector ensures that for any $j$:
\[|A_{j}|=|B_{j}|=\frac{n-1}{2}\;\;\text{and}\;\;|A_{j}\cap B_{j}|=0\]
and for $j_{1}\neq j_{2}$:
\[|A_{j_1}\cap A_{j_2}|=|A_{j_1}\cap B_{j_2}|=|B_{j_1}\cap B_{j_2}|=\frac{n-1}{4}\]
Thus, one can easily check that 
\[\mathcal A=\{A_{1}\cup\{n\},B_1\cup\{n\},\ldots, A_{n-2}\cup \{n\},B_{n-2}\cup\{n\}\}\]
is a $1$-defect $\ell$-oddtown in $[n]$ of size $2n-4$, provided $n\equiv 5 \pmod 8$ if $\ell=2$ or $\ell\mid n+3$ if $\ell>2$. Thus, a $1$-defect $\ell$-oddtown in $[n]$ of order $2n-4$ exists provided a Hadamard matrix of order $n-1$ exists and these divisibility conditions on $n$ are satisfied. We claim now that there are infinitely many values of $n$ for which this holds.

For $\ell=2$, this is ensured by a construction of Paley~\cite{P33} of Hadamard matrices of order $q+1$ for any odd prime power $q$. For $\ell>2$, this is ensured by a result of Wallis~\cite{W76} which states that for any $q\in\mathbb{N}$ there is $s_0\in \mathbb{N}$ such that a Hadamard matrix of order $2^{s}q$ exists for any $s\ge s_0$ (just take $n$ to be of the form $2^{s}q+1$, where $q=\ell-1$ and $s$ is any sufficiently large multiple of $\ell-1$). We conclude that for any prime $\ell$ there are $1$-defect $\ell$-oddtowns in $[n]$ of size $2n-4$ for infinitely many values of $n$.

\medskip

Now we prove that any $1$-defect $\ell$-oddtown in $[n]$ has size at most $\max\{n,2n-4\}$ if $\ell$ is a prime number. Suppose $\mathcal A$ is a $1$-defect $\ell$-oddtown in $[n]$. If all pairwise intersections of sets in $\mathcal A$ have size $= 0\pmod \ell$ then $\mathcal A$ is an $\ell$-oddtown and so, as discussed in the introduction, we have $|\mathcal A|\le n$. Otherwise, we can label the sets in $\mathcal A$ as $A_1,B_1,\ldots, A_{t},B_{t},A_{t+1},\ldots,A_{s}$ ($1\le t\le s$) such that the pairs $(A_i,B_i)$ have pairwise intersection of size $\neq 0\pmod \ell$ and all other pairwise intersections have size $= 0 \pmod \ell$. Let $a_1,\ldots,a_{s}$ and $b_1,\ldots, b_{t}$ be the characteristic vectors in $\mathbb F_{\ell}^{n}$ corresponding to the sets $A_1,\ldots, A_{s}$ and $B_1,\ldots, B_{t}$, respectively. Note crucially that for $i\neq j$ we have $a_i\cdot a_j = |A_i\cap A_j| = 0 \pmod \ell$, $a_i\cdot b_j = |A_i\cap B_j| = 0\pmod \ell$, $b_i\cdot b_j = |B_i\cap B_j| = 0\pmod \ell$, $a_i\cdot a_i = |A_i| \neq 0 \pmod \ell$, $a_i\cdot b_i = |A_i\cap B_i| \neq 0 \pmod \ell$ and $b_i\cdot b_i = |B_i| \neq 0 \pmod \ell$.

We consider now two separate cases:
\begin{enumerate}
\item[1)] $(a_1\cdot a_1)(b_1\cdot b_1)=(a_1\cdot b_1)^2$
\item[2)] $(a_1\cdot a_1)(b_1\cdot b_1)\neq (a_1\cdot b_1)^2$.
\end{enumerate}
In each case, we will either show that $|\mathcal A|\le n$ or we will find $s+2$ linearly independent vectors in $\mathbb{F}_{\ell}^{n}$. This then implies that if $|\mathcal A| >n$ then $s+2 \le n$ and so:
\[|\mathcal A|=s+t\le 2s\le 2(n-2)=2n-4.\]

We consider case 1 first. Let $v=(a_1\cdot a_1)b_1-(a_1\cdot b_1)a_1$ and note that $v\cdot a_i=0$ for any $i\in [s]$. Indeed, we have $v\cdot a_1 = (a_1\cdot a_1)(b_1\cdot a_1) - (a_1\cdot b_1)(a_1\cdot a_1) = 0$ and for $i>1$ we have $a_1 \cdot a_i = 0 $ and $b_1\cdot a_i = 0$, implying that $v\cdot a_i = 0$. Moreover, since $a_1$ and $b_1$ are distinct $\{0,1\}$-vectors one has $v\neq 0$, and
\[v\cdot v=(a_1\cdot a_1)\left[(a_1\cdot a_1)(b_1\cdot b_1)-(a_1\cdot b_1)^2\right]=0.\] 
Since $v\neq 0$, we can find a vector $v_{1}\in \mathbb{F}_{\ell}^{n}$ so that $v\cdot v_1\neq 0$. Define $v_{2}:=v-v_{1}$ and note that $v\cdot v_2=-v\cdot v_1$ since $v\cdot v=0$. We claim now that the vectors $a_1,\ldots, a_{s},v_{1},v_{2}$ are linearly independent. Indeed, if
\[\sum_{i=1}^{s}\alpha_{i}a_i+\beta_{1} v_{1}+\beta_{2} v_{2}=0\]
is a linear combination of these vectors, then doing the dot product with $v$ allows us to conclude that
\[0=\beta_{1}(v\cdot v_{1})+\beta_{2}(v\cdot v_2)=(\beta_1-\beta_2)(v\cdot v_1)\]
and therefore, since $v\cdot v_1 \neq 0$, we must have $\beta_1=\beta_2$. Then, since $\beta_{1}v_{1}+\beta_{2}v_{2}=\beta_1 v$, doing the dot product with $a_{i}$ for $i\in [s]$ we can deduce that $\alpha_{i}=0$. Finally, since $v\neq 0$ we can conclude then that $\beta_{1}=\beta_{2}=0$, and so the vectors $a_1,\ldots, a_s,v_1,v_2$ are linearly independent as claimed.

We now consider case 2 and assume for the moment that $t\ge 2$. We claim that the vectors $a_1,\ldots, a_{s},b_1,b_2$ are linearly independent. Indeed, if
\[\sum_{i=1}^{s}\alpha_{i}a_i+\beta_{1} b_{1}+\beta_{2} b_{2}=0\]
is a linear combination of these vectors, then doing the dot product of the above with $a_{i}\in [s]\setminus \{1,2\}$ allows us to conclude that $\alpha_{i}=0$ and so
\[\alpha_{1}a_1+\alpha_2a_2+\beta_{1}b_{1}+\beta_{2}b_2=0.\]
Now, doing the dot product of the latter with $a_1$ and $b_1$ we see that:
\[\begin{bmatrix}
(a_1\cdot a_1) & (a_1\cdot b_1)\\ 
(a_1\cdot b_1) & (b_1\cdot b_1)
\end{bmatrix}
\begin{bmatrix}
\alpha_1\\ 
\beta_1
\end{bmatrix}=\begin{bmatrix}
0\\ 
0
\end{bmatrix}\]
Since the determinant of this matrix is non-zero (because we are in case 2), we conclude that $\alpha_1=\beta_1=0$. Then, since $a_2$ and $b_2$ are distinct $\{0,1\}$-vectors we conclude that $\alpha_2=\beta_2=0$ and so the vectors $a_1,\ldots, a_s,b_1,b_2$ are linearly independent as claimed. Finally, if $t=1$ then one can show, similarly to the above, that the $s+1$ vectors $a_1,\ldots,a_{s},b_1$ are linearly independent, implying that $|\mathcal A|=s+1\le n$. This finishes the proof of Theorem~\ref{thm:1defectloddtown}.

\section{Further remarks and open problems}
\label{section:concludingremarks}

Theorems \ref{thm:vukwiseeventown}, \ref{thm:kwiseeventown} and \ref{thm:strongkwiseeventown} establish the maximum size of (strong) $k$-wise $\ell$-eventowns and characterize
 their structure for $\ell = 2$. Far less is known for (strong) $k$-wise $\ell$-eventowns with $\ell > 2$. A natural analogue of Construction~\ref{cons:maxkwiseeventown} for $\ell>2$ arises from the next construction:
\begin{cons}
\label{cons:strongkwiseleventown}
Let $B_1,\ldots, B_{\lfloor n/\ell\rfloor}$ be $\lfloor n/\ell\rfloor$ disjoint subsets of $[n]$ of size $\ell$. Then the family $\mathcal A=\left\{\bigcup_{i\in S} B_i: \;S\subseteq \left[\lfloor n/\ell\rfloor\right]\right\}$ is a strong $k$-wise $\ell$-eventown of size $2^{\lfloor n/\ell\rfloor}$ for every $k\in \NN$.
\end{cons}
Construction~\ref{cons:strongkwiseleventown} provides a strong $k$-wise $2$-eventown of maximum possible size for any $k\ge 2$, and, in light of Theorem~\ref{thm:strongkwiseeventown}, this is the unique such family for $k\ge 3$, up to the choice of the sets $B_1,\ldots, B_{\lfloor n/2\rfloor}$. Surprisingly, for $\ell > 2$, Construction~\ref{cons:strongkwiseleventown} is far from best possible. As mentioned in the introduction, Frankl and Odlyzko~\cite{FO83} constructed a strong $2$-wise $\ell$-eventown of size $2^{\Omega(\log \ell/\ell)n}$, as $n\rightarrow \infty$, which is significantly larger than the families in Construction~\ref{cons:strongkwiseleventown} for large $\ell$.

Interestingly, this phenomena does not hold only for $k=2$. Indeed, Frankl and Odlyzko's construction can be used to construct a strong $3$-wise $\ell$-eventown of size $2^{\Omega(\log \ell/\ell)n}$, as $n\rightarrow \infty$. This follows from the next simple lemma which shows how to create a large strong $k$-wise $\ell$-eventown from a large strong $(k-1)$-wise $\ell$-eventown if $k$ is odd. We leave its proof as an exercise to the interested reader.
\begin{lemma}
\label{lemma:stepup}
Suppose $\mathcal A=\{A_1,\ldots, A_m\}$ is a strong $(k-1)$-wise $\ell$-eventown on the universe $[n]$. For each $i\in [m]$ define the sets $A^{*}_{i}=(([n]\setminus A_i)+n)\subseteq [2n]\setminus [n]$ and $B_i=A_i\cup A^{*}_{i}$. If $\ell\mid n$ and $k$ is odd then $\mathcal B=\{B_1,\ldots, B_m\}$ is a strong $k$-wise $\ell$-eventown on the universe $[2n]$ of size $|\mathcal B|=|\mathcal A|$.
\end{lemma}

We can also show that for any fixed $k\in \NN$ there are strong $k$-wise $\ell$-eventowns of size $2^{\Omega(\log \ell/\ell)n}$ as $n\rightarrow \infty$ when $\ell$ is a power of $2$:

\begin{lemma}
\label{lemma:beatingblockfamilies}
For any $k\in \NN$ and $\ell$ a power of $2$, there are strong $k$-wise $\ell$-eventowns in the universe $[n]$ of size $\left(2^{k+1}\ell\right)^{\lfloor n/(2^k\ell)\rfloor}$. 
\end{lemma}
We give a brief sketch on how to construct such families. We start by recursively defining for $r\ge 0$ a family $\mathcal A^{r}$ with $2^{r+1}$ subsets $A_1^{r},\ldots, A_{2^{r+1}}^{r}$ of $[2^{r}]$ with the property that
\begin{align}
\label{eq:property}
2^{r-|S|}\; \text{  divides  }\; \left|\bigcap_{i\in S}A_i^{r}\right|  \;\text{  for any set  }\; S\subseteq [2^{r+1}]\; \text{  of size  } \;|S|\le r.
\end{align}
For $r=0$ we define $A_{1}^{0} = \emptyset$ and $A_{1}^{1}=\{1\}$. For $r> 0$, define for $1\le i\le 2^{r}$ the sets $A_{i}^{r}:= A_i^{r-1}\cup (A_{i}^{r-1}+2^{r-1})$ and $A_{i+2^{r}}^{r}=A_{i}^{r-1}\cup \left(([2^{r-1}]\setminus A_{i}^{r-1})+2^{r-1} \right)$. Finally, define $\mathcal A^{r} =\{A_1^{r},\ldots, A_{2^{r+1}}^{r}\}$. One can prove by induction on $r$ that this family satisfies property~(\ref{eq:property}).

Now, suppose $\ell = 2^{a}$. Property (\ref{eq:property}) implies that $\mathcal A^{k + a}$ is a strong  $k$-wise $\ell$-eventown in $[2^{k+a}] = [2^{k}\ell]$ of size $2^{k+a+1} = 2^{k+1}\ell$. For $j \in [\lfloor n/(2^{k}\ell)\rfloor]$ let $\mathcal B_j = \{A + (j-1)2^{k}\ell : A\in \mathcal A^{k+a}\}$ and define
\[\mathcal B = \left\{\bigcup_{j\in [\lfloor n/(2^{k}\ell)\rfloor]} B_j : B_j \in \mathcal B_j \text{ for } j\in [\lfloor n/(2^{k}\ell)\rfloor]\right\}\]
A moment's thought shows that $\mathcal B$ is a strong $k$-wise $\ell$-eventown in $[n]$ of size $|\mathcal B| = (2^{k+1}\ell)^{\lfloor n/(2^{k}\ell)\rfloor}$.

\medskip

Frankl and Odlyzko conjectured in \cite{FO83} that for any $\ell\in \NN$ there exists $k(\ell)\in \NN$ such that  if $k\ge k(\ell)$ then any $k$-wise $\ell$-eventown has size at most $2^{(1+o(1))n/\ell}$ as $n\rightarrow \infty$ (which would be asymptotically tight by Construction~\ref{cons:strongkwiseleventown}). Lemma~\ref{lemma:beatingblockfamilies} implies that if such $k(\ell)$ exists then $k(\ell) \ge (1-o(1))(\log_2\log_2\ell)$, at least when $\ell$ is a power of $2$.

Note that Lemma~\ref{lemma:beatingblockfamilies} shows that at least when $\ell$ is a power of $2$ we have strong $k$-wise $\ell$-eventowns in $[n]$ of size roughly $2^{C(k)\left(\log \ell/\ell\right) n}$, where $C(k)\sim k2^{-k}$. Moreover, if there were an analogue of Lemma~\ref{lemma:stepup} for any $k$ (not just $k$ odd) then for any $\ell \in \NN$ one could start from Frankl and Odlyzko's construction and iterate such lemma $k-2$ times in order to obtain a strong $k$-wise $\ell$-eventown in $[n]$ of size roughly $2^{C(k)\left(\log \ell/\ell\right) n}$ where $C(k) \sim 2^{-k}$. We find it plausible that such families exist for any $k,\ell\in\NN$, provided $n$ is sufficiently large (depending on $k$ and $\ell$).

In Theorem~\ref{thm:ddefectloddtown} we showed that for any $d\in \NN$ and $\ell$ a prime number, any $d$-defect $\ell$-oddtown in the universe $[n]$, for $n$ large, has size at most $(d+1)(n-2(\lceil\log_{2}(d+2)\rceil-1))$, improving  Vu's upper bound of $(d+1)n$ described at the beginning of Section~\ref{section:ddefectloddtowns}. Vu~\cite{V99} also showed that there exist $d$-defect $\ell$-oddtowns in $[n]$ of size $(d+1)(n-\ell\lceil \log_{2}(d+1)\rceil)$. These families come from the following construction:
\begin{cons}
\label{cons:ddefectloddtown}
Let $t = \lceil\log_{2}(d+1)\rceil$, $s=\ell t$ and $\mathcal S$ be a collection of $d+1$ subsets of $[t]$. Moreover, let $B_1,\ldots, B_{t}$ be $t$ disjoint subsets of $[s]$, each of size $\ell$. For each $S\in \mathcal S$ let $B_{S}=\cup_{i\in S} B_i$ and define $\mathcal B =\{B_S:S\in\mathcal S\}$. Then, the family $\mathcal A$ defined by
\[\mathcal A=\{B\cup \{i\}:B\in\mathcal B, i\in [n]\setminus [s]\}\]
is a $d$-defect $\ell$-oddtown of size $ |\mathcal A|=(d+1)(n-s)$. Indeed, for $B,B'\in\mathcal B$ and $i,i'\in [n]\setminus [s]$ we have
\[|(B\cup \{i\})\cap (B'\cup \{i'\})|=|B\cap B'|+|\{i\}\cap \{i'\}|\equiv |\{i\}\cap \{i'\}| \pmod \ell\]
and the latter is non-zero modulo $\ell$ if and only if $i= i'$.
\end{cons}
This construction can be improved for some values of $\ell$ and $d$. Notice, that the only relevant property of family $\mathcal B$ in Construction~\ref{cons:ddefectloddtown} is that it is an $\ell$-eventown on the universe $[s]$ of size at least $d+1$. Thus, if there exists an $\ell$-eventown of size $d+1$ in a universe of size smaller than $\ell\lceil \log_{2}(d+1)\rceil$ then we can improve Vu's lower bound on the maximum size of a $d$-defect $\ell$-oddtown. Frankl and Odlyzko's construction mentioned earlier shows that an $\ell$-eventown in the universe $[s]$ of size at least $2^{c\left(\log\ell/\ell\right)s}$ exists for some constant $c>0$ as $s\rightarrow \infty$. Since $2^{c\left(\log \ell /\ell\right)s}\ge d+1$ if $s\ge c^{-1}\left(\ell/\log \ell\right)\log_{2}(d+1)$, this implies that there are $d$-defect $\ell$-oddtowns of size $(d+1)(n-C\left(\ell/\log \ell\right)\log_{2}(d+1))$ for some constant $C>0$ as $n\rightarrow \infty$, provided $d$ is big enough as a function of $\ell$. It is unclear to us whether the maximum size of a $d$-defect $\ell$-oddtown should depend on $\ell$. We remark that for $d=1$, as Theorem~\ref{thm:1defectloddtown} shows, this is not the case. 

In \cite{SV05} Szab\'o and Vu considered the related problem of maximizing the size of a $k$-wise oddtown, i.e., a family of odd-sized sets such that the intersection of any $k$ has even size. They showed that if $k-1$ is a power of $2$ then for large $n$ the answer is $(k-1)(n-2\log_{2}(k-1))$. An example of a $k$-wise oddtown of this size is the one in Construction~\ref{cons:ddefectloddtown} with $d = k-2$ and $\ell = 2$. For the natural generalization of this problem modulo $\ell>2$, Szab\'o and Vu believed that Construction~\ref{cons:ddefectloddtown} with $d= k-2$ provided a $k$-wise $\ell$-oddtown in $[n]$ of maximum possible size, namely, $(k-1)(n-\ell\lceil \log_{2}(k-1)\rceil)$. This turns out not to be the case. Indeed, as described in the previous paragraph, by making a more appropriate choice of $\mathcal B$ in Construction~\ref{cons:ddefectloddtown} one can obtain for suitable values of $k$ and $\ell$ a $k$-wise $\ell$-oddtown of size $(k-1)(n-C\left(\ell/\log \ell\right)\log_{2}(k-1))$ for some constant $C>0$ and $n$ sufficiently large.

\medskip
\noindent
\textbf{Acknowledgements.} We would like to thank Shagnik Das for helpful discussions and comments.

\end{document}